\documentclass[12pt]{amsart}
\usepackage{amsmath,amssymb,amsthm}
\usepackage{mathtools}
\usepackage{stmaryrd} 
\usepackage[inline]{enumitem}
\usepackage{graphicx,setspace}
\usepackage{tikz}
\usetikzlibrary{arrows}
\usepackage{wasysym}
\usepackage[makeroom]{cancel}
\usepackage{fullpage}
\usepackage{comment}
\usepackage{enumitem}
\usepackage{caption}
\usepackage{mathtools}
\usepackage{url}
\usepackage{cleveref}
\usepackage{algorithm}
\usepackage{algorithmicx}
\usepackage{algpseudocode}
\usepackage{enumitem}

\allowdisplaybreaks

\newcommand{\ZNx}{\mathbb{Z}/N\mathbb{Z}[x]}
\newcommand{\Zpx}{\mathbb{Z}/p\mathbb{Z}[x]}
\newcommand{\Fpx}{\mathbb{F}_p[x]}

\newcommand{\Algoref}[1]{Algorithm \ref{#1}}
\newcommand{\lineref}[2]{line \algref{#1}{#2}}
\newcommand{\Lineref}[2]{Line \algref{#1}{#2}}
\newcommand{\linesref}[3]{lines \algref{#1}{#2} to \algref{#1}{#3}}
\newcommand{\Linesref}[3]{Lines \algref{#1}{#2} to \algref{#1}{#3}}

\theoremstyle{plain}
\newtheorem{theorem}{Theorem}[section]
\newtheorem{proposition}[theorem]{Proposition}

\newtheorem{lemma}[theorem]{Lemma}

\newtheorem*{claim*}{Claim}
\theoremstyle{definition}

\newtheorem{definition}[theorem]{Definition}

\newtheorem*{theorem-non}{Corollary}
\newtheorem*{definition-non}{Definition}

\def\Pr{\operatorname{Pr}}

\begin{document}

\title{On a Modification of the Agrawal-Biswas Primality Test}

\author[H. Kim]{H.\ Kim}
\email{\url{hyunjong.kim@math.wisc.edu}}
\address{Department of Mathematics \\ University of Wisconsin-Madison \\ Madison, WI 53706, USA}

\begin{abstract}
We present a variant of the Agrawal-Biswas algorithm, a Monte Carlo algorithm which tests the primality of an integer $N$ by checking whether or not $(x+a)^N$ and $x^N + a$ are equivalent in a residue ring of $\ZNx$. The variant that we present is also a randomization of Lenstra jr.~and Pomerance's improvement to the Agrawal-Kayal-Saxena deterministic primality test. We show that our variant of the Agrawal-Biswas algorithm can be used with the Miller-Rabin primality test to yield an algorithm which is slower than the Miller-Rabin test but relatively more accurate.
\end{abstract}

\maketitle

\section{Introduction}

Let $N$ be a positive integer. There are various algorithms, each with its own features, which try to determine whether or not $N$ is prime. For instance, the Fermat test, which depends on Fermat's little theorem, runs in $\tilde{O}(\log^2 N)$ time, but is prone to failure when $N$ is a Carmichael number, of which there are infinitely many \cite{AGP}. The Miller-Rabin test \cite{Miller, Rabin} modifies the Fermat test and also runs in $\tilde{O}(\log^2 N)$ time. Moreover, the Miller-Rabin test has an error probability of at most $\frac{1}{4}$ for every odd composite $N$ \cite{Monier, Rabin}. Under GRH, the Miller-Rabin test produces a deterministic primality test which runs in $\tilde{O}(\log^4 N)$ time \cite{Miller}. Elliptic curve primality proving \cite{GK} is a method with a variant \cite{AM} which heuristically runs in $\tilde{O}(\log^4 N)$ time and produces a certificate of the primality of $N$ of size $O(\log^2 N)$. This certificate can be verified in $\tilde{O}(\log^4 N)$ time. \par

The Agrawal-Kayal-Saxena [AKS] algorithm \cite{AKS} is the first primality test that was shown to be deterministic. The initial version of the algorithm has an asymptotic time complexity of $\tilde{O}(\log^{12} N)$, where $N$ is the positive integer whose primality is tested. The authors proposed another version of the algorithm which has an asymptotic time complexity of $\tilde{O}(\log^{10.5} N)$. Later, \cite{LP} demonstrated yet another variant of the AKS algorithm that runs in $\tilde{O}(\log^6 N)$ time using pseudofields. \par

We modify the Agrawal-Biswas primality test, upon which the AKS algorithm is based. The resulting Monte Carlo algorithm can be combined with the Miller-Rabin algorithm to yield an algorithm that runs in $\tilde{O}(\log^{c+2} N)$ time, for any sufficiently large $N$ and any $c > \frac{46}{25}$, and fails for composite $N$ with a probability less than $\frac{1}{2^{\Omega(\log^{c+1} N)}}$. \par

\Cref{section:AB} summarizes the Agrawal-Biswas algorithm. \Cref{section:acc} and \Cref{section:run} respectively analyze the accuracy and runtime of the modified Agrawal-Biswas algorithm that we provide. We then discuss the Miller-Rabin algorithm in \Cref{section:mr} to understand how the Miller-Rabin algorithm and the modified Agrawal-Biswas combine in \Cref{section:both}. For a summary of the runtime and accuracy of this combined algorithm, see \Cref{theorem:summary}. In \Cref{section:compare}, we compare the accuracy guarantees of the Miller-Rabin algorithm and the combined algorithm from \Cref{section:both} relative to their runtimes. We use details from \cite{LP} to discuss pseudofields in \Cref{section:polyconstruct} and to show that pseudofields constructed by an algorithm in \cite{LP} are fields when reduced modulo prime factors of $N$. In particular, these pseudofields can be used to obtain polynomials which yield the high accuracy of the combined algorithm guaranteed in \Cref{proposition:strongaccuracy}.\par

For more on runtime analysis of algorithms for basic integer and polynomial arithmetic, see \cite{ModernComputerAlgebra}.

\section{The Agrawal-Biswas Algorithm} \label{section:AB}

The Agrawal-Biswas algorithm is based on the following lemma:
\begin{lemma}{\cite[Lemma 3.1]{AB}} \label{lemma:main}
	An integer $N > 1$ is prime if and only if, for any $a \in \mathbb{Z}/N\mathbb{Z}^\times$, the equality 
	\begin{align*}
		(x+a)^N = x^N + a
	\end{align*}
	holds in the polynomial ring $\mathbb{Z}/N\mathbb{Z}[x]$. 
	\begin{proof}
		Let $a \in \mathbb{Z}/N\mathbb{Z}^\times$. By the binomial theorem, 
		\begin{align*}
			(x+a)^N = \sum_{i=0}^N {n \choose i} x^{N-i} a^i.
		\end{align*}
		Suppose that $N$ is prime. For all $i = 1,\ldots,N-1$, the binomial coefficient ${N \choose i}$ is divisible by $N$ and is hence $0$ in $\mathbb{Z}/N\mathbb{Z}$. Therefore, 
		\begin{align*}
			(x+a)^N = x^N + a^N. 
		\end{align*}
		By Fermat's little theorem, $a^N = a$, so $(x+a)^N = x^N + a$ as desired. \par
		Suppose that $N$ is instead composite. Let $p$ be any prime factor of $N$ and say that $p^\alpha \mid \mid N$. The binomial coefficient ${n \choose n-p}$ is not divisible by $p^\alpha$. Moreover, $a$ is assumed to be in $\mathbb{Z}/N\mathbb{Z}^\times$, so $(x+a)^N$ has a nonzero $x^{N-p}$ term whose coefficient is ${N \choose N-p} a^p$. Thus, $(x+a)^N \neq x^N+a$. 
	\end{proof}
\end{lemma}

Thus, testing the primality of $N$ is equivalent to $(x+A)^N = x^N+a$. However, testing $(x+a)^N = x^N + a$ by computing the power $(x+a)^N$ is computationally infeasible. One can instead test $(x+a)^N=x^N+a$ modulo a polynomial $f(x) \in \mathbb{Z}/N\mathbb{Z}[x]$, deeming $N$ to be prime if $(x+a)^N \equiv x^N+a \pmod{f(x)}$ and composite otherwise. Although this test is computationally feasible and completely accurate when $N$ is prime, it is also prone to error if $N$ is composite. \par

The algorithm in \cite[Section 3]{AB} chooses $f$ randomly of degree $\lceil \log N \rceil$. If $N$ is a composite number that is not a perfect power and whose prime factors are greater than $13$, then \cite[Theorem 3.2]{AB} guarantees $\frac{1}{3}$ as an upper bound to the probability that the test fails. Difficulties in asymptotically improving this bound arise because $f$ might be a reducible polynomial modulo a prime factor $p$ of $N$. \par

\subsection{Modifying the Agrawal-Biswas Algorithm}

Instead of letting $f(x) \in \ZNx$ be randomly chosen, we let $f(x)$ be a polynomial of sufficiently large degree and irreducible modulo prime factors $p$ of $N$. Moreover, instead of checking $(x+1)^N \equiv x^N+1 \pmod{f(x)}$, we check that $(h(x) + 1)^N \equiv h(x)^N + 1 \pmod{f(x)}$ for a randomly chosen $h(x) \in \ZNx$ of degree less than $\deg f$. We also check if $N$ has small prime factors to decrease the probability that the test fails, see \Cref{proposition:accuracy2}. 

\begin{algorithm} \caption{Hi}
\begin{algorithmic}[1]
	\Function{PrimalityTestAB}{$N$,$f(x)$}
		\Require{$N > 1$, $f(x)$ is irreducible modulo some prime factor $p$ of $N$}
		\For{$2 <= a <= \log N$} \label{line:smallprimestart}
			\If{$a$ divides $N$}
				\State \Return{COMPOSITE}
			\EndIf
		\EndFor \label{line:smallprimeend}
		\State $h(x) \gets \text{uniformly random polynomial in $\ZNx$ of degree at most $\deg f - 1$}$ \label{line:random}
		\If{$(h(x)+1)^N \equiv h\left( x \right)^N + 1 \pmod{f(x)}$} \label{line:compute}
			\State \Return{PRIME}
		\Else
			\State \Return{COMPOSITE}
		\EndIf
	\EndFunction
\end{algorithmic} \label{algorithm:AB}
\end{algorithm}


\section{Accuracy Analysis of \Algoref{algorithm:AB}} \label{section:acc}

If $N$ is prime, then \Algoref{algorithm:AB} always correctly determines $N$ to be prime by \Cref{lemma:main}. On the other hand, \Algoref{algorithm:AB} might incorrectly determine $N$ to be prime when $N$ is composite. \Cref{proposition:prop} applies the Schwartz-Zippel lemma to bound the probability that \Algoref{algorithm:AB} fails. 

\begin{proposition} \label{proposition:prop}
	Let $p$ be a prime, let $f(x) \in \Zpx$ be irreducible, and let $g(x) \in \Zpx$ be nonzero. Let $h(x) \in \Zpx$ be a polynomial chosen uniformly at random such that $\deg h < \deg f$. Then, 
	\begin{align*}
		\Pr_{h(x)}[g(h(x)) \equiv 0 \pmod{f(x)} ] \leq \frac{\deg g}{p^{\deg f}}.
	\end{align*}
	\begin{proof}
		Since $f$ is irreducible, the residue ring $\Fpx/(f(x))$ is isomorphic to $\mathbb{F}_{p^{\deg f}}$, the finite field of $p^{\deg f}$ elements. In particular, the image $\overline{x}$ of $x$ in $\Fpx/(f(x))$ corresponds to a generator $\alpha$ of $\mathbb{F}_{p^{\deg f}}$ whose minimal polynomial over $\mathbb{Z}/p\mathbb{Z} \simeq \mathbb{F}_p$ is $f$. Since $\deg h < \deg f$, the image $h(\overline{x})$ of $h(x)$ in $\Fpx/(f(x))$ corresponds to $h(\alpha)$ in $\mathbb{F}_{p^{\deg f}}$, so choosing $h$ uniformly at random corresponds to choosing an element of $\mathbb{F}_{p^{\deg f}}$ uniformly at random. Furthermore, the event that $g(h(x)) \equiv 0 \pmod{f(x)}$ if equivalent to the event that $g(h(\alpha)) = 0$, i.e. that $h(\alpha)$ is a root of $g$. Since $g$ has at most $\deg g$ roots and $\mathbb{F}_{p^{\deg f}}$ has $p^{\deg f}$ elements,  
		\begin{align*}
			\Pr_{h(x)}[g(h(x)) \equiv 0 \pmod{f(x)} ] \leq \frac{\deg g}{p^{\deg f}} = \Pr_{\beta \in \mathbb{F}_{p^{\deg f}}}[g(\beta) = 0] \leq \frac{\deg g}{p^{\deg f}}. 
		\end{align*}
	\end{proof}
\end{proposition}

\begin{proposition} \label{proposition:accuracy}
	If $N$ is composite, then \Algoref{algorithm:AB} fails with probability less than  
	\begin{align*}
		\frac{N}{p^{\deg f}}
	\end{align*}
	where $p$ is any prime factor of $N$ for which $f(x)$ is irreducible modulo $p$. In particular, if $N$ has no prime factors less than $\log N$ as well, then \Algoref{algorithm:AB} fails with probability less than
	\begin{align*}
		\frac{N}{\log^{\deg f} N}.
	\end{align*}
	\begin{proof}
		\Lineref{algorithm:AB}{line:compute} tests $g(h(x)) \equiv 0 \pmod{f(x)}$ where $g(x) = (x+1)^N - (x^N + 1)$. Note that $g(x)$ is a polynomial of degree less than $N$. \Algoref{algorithm:AB} fails exactly when the equivalence $g(h(x)) \equiv 0 \pmod{f(x)}$ holds over $\mathbb{Z}/N\mathbb{Z}$, in which case the equivalence holds over $\mathbb{Z}/p\mathbb{Z}$. Since $h(x)$ is chosen to be a polynomial over $\mathbb{Z}/N\mathbb{Z}$ uniformly at random, its residue modulo $p$ is a uniformly random polynomial over $\mathbb{Z}/p\mathbb{Z}$. Therefore, by \Cref{proposition:prop}, \Algoref{algorithm:AB} fails with probability at most
		\begin{align*}
			\frac{\deg g}{p^{\deg f}}
		\end{align*}
		which is less than
		\begin{align*}
			\frac{N}{p^{\deg f}}.
		\end{align*}
		If $N$ has no prime factors less than $\log N$ as well, then this is less than
		\begin{align*}
			\frac{N}{\log^{\deg f} N}.
		\end{align*}
	\end{proof}
\end{proposition}

If $f(x)$ is irreducible modulo all prime factors of $N$, then we can strengthen the accuracy guarantee of \Algoref{algorithm:AB}. 
\begin{proposition} \label{proposition:accuracy2}
	Suppose that $N$ is a composite number whose prime factorization is $N = \prod_{i=1}^r p_i^{e_i}$ where the $p_i$ are distinct prime numbers and $e_i \geq 0$. The probability that \Algoref{algorithm:AB} fails is less than
	\begin{align*}
		\frac{N^r}{\prod_{i=1}^r p_i^{\deg f}}.
	\end{align*}
	\begin{proof}
		Just as in \Cref{proposition:accuracy}, \Algoref{algorithm:AB} fails exactly when the equivalence $g(h(x)) \equiv 0 \pmod{f(x)}$ holds over $\mathbb{Z}/N\mathbb{Z}$ where $g(x) = (x+1)^N - (x^N+1)$. In particular, the equivalence holds over $\mathbb{Z}/p_i\mathbb{Z}$ for every $i$. Since $h(x)$ is chosen to be a polynomial over $\mathbb{Z}/N\mathbb{Z}$ of degree less than $\deg f$ uniformly at random, the tuple
		\begin{align*}
			(h(x) \pmod{p_i})_{i=1}^n 
		\end{align*}
		of reductions of $h(x)$ modulo $p_i$ takes as values all tuples of polynomials over $\mathbb{Z}/p_i\mathbb{Z}$ of degree less than $\deg f$ with uniform probability. Just as in \Cref{proposition:accuracy}, the probability that the equivalence $g(h(x)) \equiv 0 \pmod{f(x)}$ holds over $\mathbb{Z}/p_i\mathbb{Z}$ is less than
		\begin{align*}
			\frac{N}{p^{\deg f}}.
		\end{align*}
		Therefore, the probability that the equivalence holds over $\mathbb{Z}/p_i\mathbb{Z}$ for every $i$ is less than
		\begin{align*}
			\frac{N^r}{p^{\deg f}}.
		\end{align*}
	\end{proof}
\end{proposition}

\section{Runtime Analysis of \Algoref{algorithm:AB}} \label{section:run}

We analyze the runtime of \Algoref{algorithm:AB}. 
\begin{proposition} \label{proposition:runtime1}
	In \Algoref{algorithm:AB}, let $D = \deg f$. Assuming that each random bit can be generated in $O(1)$ time, \Algoref{algorithm:AB} runs in $\tilde{O} \left( D \log^2 N  \right)$ time and requires $O(D \log N)$ random bits. In particular, if $D = \tilde{O}(\log^c N)$, then \Algoref{algorithm:AB} runs in $\tilde{O}\left( \log^{c+2} N \right)$ time and requires $O(\log^{c+1} N)$ random bits. \par
		\begin{proof}
			\Linesref{algorithm:AB}{line:smallprimestart}{line:smallprimeend} run in $\tilde{O}(\log^2 N)$ time since checking whether or not $a$ divides $N$ takes $\tilde{O}(\log N)$ time for each $a$. Moreover, \lineref{algorithm:AB}{line:random} runs in $O(D \log N)$ time and uses $O(D \log N)$ random bits. Computing $(h(x)+1)^N$ and $h(x)^N$ for \Lineref{algorithm:AB}{line:compute} can be done via $O(\log N)$ multiplications of elements in $\ZNx/(f(x))$ via binary exponentiation. Each such multiplication takes $\tilde{O}(D\log N)$ time, so computing the powers takes $\tilde{O}(D\log^2 N)$ time. Therefore, $\tilde{O} \left( D \log^2 N \right)$ time and $O(D \log N)$ random bits are needed to run \Algoref{algorithm:AB}
	\end{proof}
\end{proposition}

Recall that \Algoref{algorithm:AB} requires a polynomial $f(x) \in \ZNx$ that is irreducible modulo some prime factor of $N$. We will later observe in \Cref{section:polyconstruct} that a deterministic algorithm by \cite{LP} either correctly declares $N$ to be composite or produces a polynomial $f(x) \in \ZNx$ that is irreducible modulo all prime factors of $N$ \footnote{The algorithm only works for sufficiently large values of $N$ and produces a polynomial $f(x)$ of degree $D = \Omega\left( \log^{\frac{46}{25}} N \right)$, see \Cref{proposition:polyconstructruntime}.}. Moreover, \Cref{proposition:polyconstructruntime} shows that this algorithm runs in $\tilde{O}(D \log N)$ time, which is less time than the time needed to run \Algoref{algorithm:AB}. 

\section{The Miller-Rabin Algorithm} \label{section:mr}

The Miller-Rabin primality test is a Monte Carlo algorithm that always correctly determines prime numbers to be prime and may incorrectly determine composite numbers to be prime, just as the Agrawal-Biswas test does. Let $s$ and $t$ be nonnegative integers such that $N-1 = 2^s t$ with $t$ odd. For a nonzero element $a$ of $\mathbb{Z}/N\mathbb{Z}$, the Miller-Rabin test checks whether or not  
\begin{enumerate}
	\item $a^t \equiv 1 \pmod{N}$ \text{ or }
	\item $a^{2^i t} \equiv -1 \pmod{N}$ for some integer $i$ where $0 \leq i \leq s-1$. 
\end{enumerate}
In particular, if $N$ is prime, then at least one of these two conditions must hold. If neither of these two conditions holds, then $a$ is said to be a \textit{witness} for the compositeness of $N$. \par

\Algoref{algorithm:millerrabin} below describes pseudocode for the Miller-Rabin test.

\begin{algorithm}\caption{Miller-Rabin Primality Test}
\begin{algorithmic}[1]
	\Function{PrimalityTestMR}{$N$}
	\Require{$N > 1$}
	\If{$N == 2$} 
		\State \Return{PRIME}
	\ElsIf{$N \equiv 0 \pmod{2}$}
		\State \Return{COMPOSITE} 
	\EndIf \label{line:eventest}
	\State $a \gets \text{uniformly random integer in $[1,n-1]$}$
	\State $s,t \gets \text{nonnegative integers such that } N-1 = 2^st \text{ with } t \text{ odd}$
	\If{$a^t \equiv 1 \pmod{N}$} \label{line:powerstart}
		\State \Return{PRIME}
	\EndIf	
	\For{$1 \leq i \leq s-1$} 
		\If{$a^{2^i t} \equiv -1 \pmod{N}$}
			\State \Return{PRIME}
		\EndIf
	\EndFor \label{line:powerend}
	\State \Return{COMPOSITE}
	\EndFunction
\end{algorithmic} \label{algorithm:millerrabin}
\end{algorithm}

\subsection{Accuracy of the Miller-Rabin Primality Test} \label{subsection:mraccurate}

For a fixed odd composite number $N$, the probability that the Miller-Rabin primality test fails is bounded above by $\frac{1}{4}$. \cite[Theorem 2]{DLP} suggests that this bound is weak for general $N$ --- if $N$ is a uniformly random odd integer in $[2^{k-1},2^k]$, then the probability that $N$ is composite given that $a$ is not a witness is bounded above by
\begin{align*}
	k^2 \cdot 4^{2-\sqrt{k}}.
\end{align*}

We prove \Cref{lemma:mraccuracy} to bound the probability of error for the Miller-Rabin test given the prime factorization of $N$. Although \Cref{lemma:mraccuracy} in itself does not yield the $\frac{1}{4}$ accuracy bound for semiprime $N$, we will still use \Cref{lemma:mraccuracy} later in \Cref{proposition:strongaccuracy} to prove an accuracy bound for \Algoref{algorithm:both}, which uses both the Miller-Rabin test and the modified Agrawal-Biswas test.
\begin{lemma} \label{lemma:mraccuracy}
	Suppose that $N$ is an odd composite number whose prime factorization is $N = \prod_{i=1}^r p_i^{e_i}$ where the $p_i$ are distinct prime numbers and $e_i \geq 0$. The probability that the Miller-Rabin test fails is at most
	\begin{align*}
		\frac{1}{2^{r-1} \cdot \prod_{i=1}^r p_i^{e_i-1}}.
	\end{align*}
	\begin{proof}
		The Miller-Rabin test always succeeds when the randomly chosen nonzero base $a \in \mathbb{Z}/N\mathbb{Z}$ shares factors with $N$. Assume that $a \in \mathbb{Z}/N\mathbb{Z}^\times$. For a positive integer $m$, let $C_m$ denote the cyclic group of order $m$. Note that the multiplicative group $\mathbb{Z}/N\mathbb{Z}^\times$ is isomorphic to 
		\begin{align*}
			C := \prod_{i=1}^r C_{p_i^{e_i-1}} \times C_{p_i-1} 
		\end{align*}
		and say that $a$ corresponds to the element of $C$ whose component in $C_{p_i^{e_i-1}} \times C_{p_i-1}$ is $(x_i, y_i)$ where $x_i \in C_{p_i^{e_i-1}}$ and $y_i \in C_{p_i-1}$. If $a$ is a nonwitness, then we must have $a^{N-1} \equiv 1 \pmod{N}$, so $(N-1)x_i$ is the identity $0$ of $C_{p_i^{e_i-1}}$. Moreover, $p_i^{e_i-1}$ and $N-1$ share no common factors, so $x_i$ must be the identity if $a$ is a nonwitness. \par
		Write $N-1 = 2^s t$ where $s$ and $t$ are nonnegative integers with $t$ odd. If $a$ is a nonwitness, then
		\begin{enumerate}
			\item $a^t \equiv 1 \pmod{N}$ or 
			\item $a^{2^j t} \equiv -1 \pmod{N}$ for some integer $j$ where $0 \leq j \leq s-1$. 
		\end{enumerate}
		Furthermore, since $x_i$ is $0$, these conditions are equivalent to
		\begin{enumerate}
			\item $ty_i = 0$ for every $i$ or
			\item $2^j t y_i$ has order $2$ for every $i$ for some integer $j$ where $0 \leq j \leq s-1$.
		\end{enumerate}
		The second of these conditions is equivalent to $ty_i$ having the same order $2^{j+1}$ for every $i$. In particular, $ty_i$ is an element of the $2$-Sylow subgroup of $C$. The $2$-Sylow subgroup of $C$ is isomorphic to
		\begin{align*}
			\prod_{i=1}^r C_{|p_i-1|_2^{-1}},
		\end{align*}
		where $| p_i-1 |_2^{-1}$ is the largest power of $2$ dividing $p_i-1$. For each $i$, at most half of the elements of $C_{|p_i-1|_2^{-1}}$ can have any particular order. For a randomly chosen $a \in \mathbb{Z}/N\mathbb{Z}^\times$, the probability that $ty_i$ has the same order $2^{j+1}$ for every $i$ is at most $\frac{1}{2^{r-1}}$; if $ty_1$ has some order $2^{j+1}$, then each other $ty_i$ has that same order with probability at most $\frac{1}{2}$ because $t$ is odd. Therefore, the probability that a randomly chosen element $a$ of $\mathbb{Z}/N\mathbb{Z}^\times$ is a nonwitness for $N$ is at most 
		\begin{align*}
			\frac{1}{2^{r-1} \cdot \prod_{i=1}^r p_i^{e_i-1}}
		\end{align*}
		as desired. 
	\end{proof}
\end{lemma}

We heuristically argue that a constant bound may be the tightest possible upper bound for the probability of error for general fixed odd composite numbers $N$. We make this argument for one class of composite numbers $N$. The argument can be generalized to other classes of composite numbers. \par

Assume that there are infinitely many odd positive integers $k$ such that $p = 2k+1$ and $q = 6k+1$ are both prime. Let $N = pq$ so that $N-1 = (2k+1)(6k+1) = 12k^2+8k+1$. Note that $N-1 = 4k(3k+2)$, so expressing $N-1 = 2^s t$, we have that $k$ divides $t$. The multiplicative group $\mathbb{Z}/N\mathbb{Z}^\times$ is isomorphic to 
$$\mathbb{Z}/p\mathbb{Z}^\times \times \mathbb{Z}/q\mathbb{Z}^\times \simeq \mathbb{Z}/(p-1)\mathbb{Z} \oplus \mathbb{Z}/(q-1) \mathbb{Z} \simeq \mathbb{Z}/(2k)\mathbb{Z} \oplus \mathbb{Z}/(6k)\mathbb{Z}.$$
Since $k$ divides $t$, the elements of $\mathbb{Z}/(2k)\mathbb{Z} \oplus \mathbb{Z}/(6k)\mathbb{Z}$ of the form $(2\alpha, 6\beta)$ where $\alpha \in \mathbb{Z}/(2k)\mathbb{Z}$ and $\beta \in \mathbb{Z}/(6k)\mathbb{Z}$ become the identity element when multiplied by $t$. Thus, at least $\frac{1}{12}$ of the elements of $\mathbb{Z}/(2k)\mathbb{Z} \oplus \mathbb{Z}/(6k)\mathbb{Z}$ become the identity when multiplied by $t$. These elements correspond to elements $a$ of $\mathbb{Z}/N\mathbb{Z}^\times$ such that $a^t \equiv 1 \pmod{N}$. Thus, at least $\frac{1}{12}$ of the elements of $\mathbb{Z}/N\mathbb{Z}^\times$ are nonwitnesses. Since $\left (1 - \frac{1}{p} \right) \left( 1 - \frac{1}{q} \right)$ of the elements of $\mathbb{Z}/N\mathbb{Z}$ are elements of $\mathbb{Z}/N\mathbb{Z}^\times$, the probability that the Miller-Rabin test fails is at least
\begin{align*}
	\frac{1}{12} \left( 1 - \frac{1}{p} \right) \left( 1 - \frac{1}{q} \right)
\end{align*}
 for this class of $N$. Note that $p$ is at least $3$ and $q$ is at least $7$, so the test fails with a probability of at least
\begin{align*}
	\frac{1}{12} \cdot \frac{2}{3} \cdot \frac{6}{7} = \frac{1}{21}. 
\end{align*}

\subsection{Runtime Analysis of the Miller-Rabin Primality Test} \label{subsection:mrruntime}

\begin{proposition} \label{proposition:mrruntime}
Assuming that each random bit can be generated in $O(1)$ time, the Miller-Rabin primality test runs in $\tilde{O}(\log^2 N)$ time and requires $O(\log N)$ random bits. 
\begin{proof}
	The runtime of \Algoref{algorithm:millerrabin} is dominated by \linesref{algorithm:millerrabin}{line:powerstart}{line:powerend}. If we compute $2^{N-1}$ via left-to-right binary exponentation, then we compute $2^t$ and $2^{2^i t}$ for $1 \leq i \leq s-1$ in intermediary steps. Computing $2^{N-1}$ via binary exponentiation takes $O(\log N)$ multiplications and each multiplication takes $\tilde{O}(\log N)$ time. Thus, the algorithm runs in $\tilde{O}(\log^2 N)$ time. Moreover, all of the randomness of the primality test comes from generating a base $a \in \mathbb{Z}/N\mathbb{Z}$, which takes $O(\log N)$ random bits. 
\end{proof}
\end{proposition}

The Miller-Rabin test may terminate upon computing $a^{2^i t}$ where $0 \leq i \leq s-1$ and $2^i t$ is much smaller in comparison to $N-1$. However, this is not always the case and in fact, there may be classes of composite odd numbers $N-1$ such that the probably that the Miller-Rabin test fails is at least a constant and such that $t$ is asymptotically comparable to $N-1$. For instance, as presented in \Cref{subsection:mraccurate}, let $N = pq$ where $p$ and $q$ are prime of the form $p = 2k+1$ and $q = 6k+1$ for odd $k$. In this case, $N-1 = 4k(3k+2)$, and since $k$ is odd, $3k+2$ is odd, so $t = k(3k+2) = \frac{N-1}{4}$.

\section{Using the Modified Agrawal-Biswas and Miller-Rabin Tests Together} \label{section:both}

Let $N$ be odd and composite. The Miller-Rabin test's accuracy increases when $N$ has more prime factors, especially repeated ones, by \Cref{lemma:mraccuracy}. On the other hand, the modified Agrawal-Biswas test's accuracy increases when $N$ has larger prime factors by \Cref{proposition:accuracy2}. We can use these two tests together to obtain an algorithm with stronger probability bounds than the probability bounds guaranteed for either test. \par
In particular, we will construct this algorithm by running the Miller-Rabin test multiple times and the modified Agrawal-Biswas test just once so that the time used to run the Miller-Rabin tests and time used to run the single modified Agrawal-Biswas test are asymptotically similar. This will keep the time for the entire test asymptotically minimal while maximizing the test's accuracy relative to its runtime. Since each Miller-Rabin test runs in $\tilde{O}(\log^2 N)$ time by \Cref{proposition:mrruntime} and the modified Agrawal-Biswas runs in $\tilde{O}\left(\log^{c+2} N \right)$ time by \Cref{proposition:runtime1}, we will invoke the Miller-Rabin test $\Theta(\log^{c} N)$ times. \par

We specify the details of our proposed scheme of using the Miller-Rabin test and the Agrawal-Biswas test together in \Algoref{algorithm:both} below. To summarize, \Algoref{algorithm:both} tests the Miller-Rabin test $\deg f$ times and the modified Agrawal-Biswas test once and correctly determines $N$ to be composite if any of these tests indicate $N$ to be composite and, possibly incorrectly, indicates $N$ to be prime otherwise.
\begin{algorithm} \caption{}
\begin{algorithmic}[1]
	\Function{PrimalityTestBoth}{$N$,$f(x)$}
		\Require{$N > 1$, $f(x)$ is irreducible modulo some prime factor $p$ of $N$}
		\For{$1 \leq i \leq \deg f$} \label{line:bothMRstart}
			\State $MRResult \gets \textproc{PrimalityTestMR}(N)$
			\If{$MRResult == \text{COMPOSITE}$}
				\State \Return{COMPOSITE}
			\EndIf
		\EndFor \label{line:bothMRend}
		\State \Return{$\textproc{PrimalityTestAB}(N,f(x))$} \label{line:bothAB}
	\EndFunction
\end{algorithmic} \label{algorithm:both}
\end{algorithm}

We analyze the accuracy of \Algoref{algorithm:both} in \Cref{proposition:strongaccuracy} below.
\begin{proposition} \label{proposition:strongaccuracy}
	Suppose that $N$ is an odd composite number with $r$ distinct prime factors and that $f(x)$ is irreducible modulo $p$ for every prime factor $p$ of $N$. Further assume that the randomness of the calls to \Algoref{algorithm:AB} (the modified Agrawal-Biswas test) and \Algoref{algorithm:millerrabin} (the Miller-Rabin test) are independent. \Algoref{algorithm:both} fails with probability less than 
	\begin{align*}
		\frac{1}{2^{(r-1)\deg f} \cdot N^{\deg f - r}}.
	\end{align*}
\begin{proof}
	Let the prime factorization of $N$ be
	\begin{align*}
		N = \prod_{i=1}^{r} p_i^{e_i}
	\end{align*}
	where the $p_i$ are distinct primes and the $e_i$ are positive integers. By \Cref{lemma:mraccuracy}, each Miller-Rabin test fails with probability at most
	\begin{align*}
		\frac{1}{2^{r-1} \cdot \prod_{i=1}^r p_i^{e_i-1}}.
	\end{align*}
	Since $f(x)$ is irreducible modulo every prime factor of $N$, the modified Agrawal-Biswas test fails with probability at less than
	\begin{align*}
		\frac{N^r}{\prod_{i=1}^r p_i^{\deg f}}
	\end{align*}
	by \Cref{proposition:accuracy2}. Moreover, the randomness in the calls of the Miller-Rabin and modified Agrawal-Biswas tests are assumed to be independent, so \Algoref{algorithm:both} fails with probability less than
	\begin{align*}
		&\left( \frac{1}{2^{r-1} \cdot \prod_{i=1}^r p_i^{e_i-1}} \right)^{\deg f} \cdot \frac{N^r}{\prod_{i=1}^r p_i^{\deg f}} \\
		&= \frac{N^r}{2^{(r-1) \deg f} \prod_{i=1}^r p_i^{e_i \deg f}} \\
		&= \frac{1}{2^{(r-1) \deg f} \cdot N^{\deg f -r}}.
	\end{align*}
\end{proof}
\end{proposition}

We further analyze the runtime of \Algoref{algorithm:both} in \Cref{proposition:strongruntime} below.
\begin{proposition}\label{proposition:strongruntime}
	In \Algoref{algorithm:both}, let $D = \deg f$. Assuming that each random bit can be generated in $O(1)$ time, \Algoref{algorithm:AB} runs in $\tilde{O} \left( D \log^2 N  \right)$ time and requires $O(D \log N)$ random bits. In particular, if $D = \tilde{O}(\log^c N)$, then \Algoref{algorithm:AB} runs in $\tilde{O}\left( \log^{c+2} N \right)$ time and requires $O(\log^{c+1} N)$ random bits. \par
	\begin{proof}
		Since each Miller-Rabin test runs in $\tilde{O}(\log^2 N)$ time by \Cref{proposition:mrruntime}, \linesref{algorithm:both}{line:bothMRstart}{line:bothMRend} run in $\tilde{O}\left( D \log^2 N \right)$ time. Moreover, by \Cref{proposition:runtime1}, \lineref{algorithm:both}{line:bothAB} runs in $\tilde{O} \left( D \log^2 N \right)$ time. Therefore, \Algoref{algorithm:both} runs in $\tilde{O} \left( D \log^2 N \right)$ time. \par
		Furthermore, \linesref{algorithm:both}{line:bothMRstart}{line:bothMRend} requires $O(D \log N)$ random bits and \lineref{algorithm:both}{line:bothAB} also requires $O(D \log N)$ random bits. Therefore, \Algoref{algorithm:both} requires $O(D \log N)$ random bits.
	\end{proof}
\end{proposition}

\section{Comparing the Runtime-Accuracy Payoffs for the Modified Agrawal-Biswas and Miller-Rabin Tests} \label{section:compare}
Suppose that a Monte Carlo algorithm with input $N$ runs in at most $T(N)$ time and fails with probability at most $\epsilon(N)$. One can run the algorithm multiple times to increase the probability that the algorithm succeeds at least once. This is useful for the Agrawal-Biswas and Miller-Rabin primality tests because $N$ must be composite if the tests indicate that $N$ is composite even once. Assuming that each invocation of the algorithm is independent of the others, running the algorithm $t$ times extends the algorithm to one which runs in at most $tT(N)$ time and fails (for each of the $t$ invocations) with probability at most $\epsilon(N)^t$. \par

A slower probabilistic algorithm, which we call algorithm $1$, can be advantageous over a faster one, which we call algorithm $2$, if a sufficiently high degree of confidence is desired and each invocation of algorithm $1$ is sufficiently more accurate than each invocation of algorithm $2$. Suppose that algorithm $i$, for $i = 1,2$ has runtime bounded above by $T_i(N)$ and failure probability bounded above by $\epsilon_i(N)$. Further say that we want to run each algorithm enough times to ensure that the probability of failure is less than $\delta$. In this case, we need to invoke algorithm $i$ at least
\begin{align*}
	\left \lceil \frac{\log \delta}{\log \epsilon_i(N)} \right \rceil
\end{align*}
times\footnote{We assume that the base of the logarithms is $2$ and, in particular, is greater than $1$. Thus, $\log \delta$ and $\log \epsilon_i(N)$ are both negative.}. The total time used in invoking algorithm $i$ this many times is  
\begin{align*}
	O\left( \frac{\log \delta}{\log \epsilon_i(N)} \cdot T_i(N) \right).
\end{align*}
Given that $T_i(N)$ and $\epsilon_i(N)$ are tight bounds, in the sense that $T_i(N)$ and $\log \epsilon_i(N)$ are accurate within constant factors, and that
\begin{align*}
	\frac{T_1(N)}{|\log \epsilon_1(N)|} = \omega\left(  \frac{T_2(N)}{|\log \epsilon_2(N)|} \right),
\end{align*}
i.e.
\begin{align*}
\frac{T_1(N)}{|\log \epsilon_1(N)|}
\end{align*}
is asymptotically greater than
\begin{align*}
\frac{T_2(N)}{|\log \epsilon_2(N)|},
\end{align*}
it takes asymptotically less time to achieve the desired accuracy $\delta$ by repeating algorithm $2$ than by repeating algorithm $1$. In this sense, a probabilistic algorithm with runtime $T(N)$ and failure probability $\epsilon(N)$ can be considered to be asymptotically accurate relative to the runtime when the fraction
\begin{align*}
	\frac{T(N)}{|\log \epsilon(N)|}
\end{align*}
is small.

Let $T_{AB}(N)$ and $T_{MR}(N)$ respectively be the (upper bounds for the) runtimes of \Algoref{algorithm:both} and the Miller-Rabin primality test. Similarly, let $\epsilon_{AB}(N)$ and $\epsilon_{MR}(N)$ respectively be (upper bounds for the) probabilities that \Algoref{algorithm:both} and the Miller-Rabin test fail when $N$ is composite. In \Algoref{algorithm:AB} and \Algoref{algorithm:both}, assume that $f(x) \in \ZNx$ is irreducible over $\mathbb{Z}/p\mathbb{Z}$ for every prime factor $p$ of $N$ and specify $D = \deg f = \Theta(\log^c N)$ for a constant $c > \frac{46}{25}$. Furthermore, let $r$ be the number of distinct prime factors of $N$, so $r = O(\log N)$. By \Cref{proposition:runtime1} and \Cref{proposition:strongaccuracy}, 
\begin{align*}
	T_{AB}(N) = \tilde{O}\left( \log^{c+2} N \right)
\end{align*}
and
\begin{align*} 
	\epsilon_{AB}(N) < \frac{1}{2^{(r-1) \log^c N} N^{\log^c N - r}}.
\end{align*}
Moreover,
\begin{align*}
	T_{MR}(N) = \tilde{O}(\log^{2} N )
\end{align*}
and
\begin{align*}
	\epsilon_{AB}(N) < \frac{1}{4}.
\end{align*}
Comparing the asymptotic accuracies of the two algorithms relative to their runtimes is thus tantamount to comparing the fractions
\begin{align*}
	R_{AB}(N) := \frac{T_{AB}(N)}{|\log \epsilon_{AB}(N)|} = \tilde{O}\left( \frac{\log^{c+2}N }{(r-1) \log^c N + (\log^c N-r)\log N} \right) = \tilde{O}\left( \frac{\log^{c+2} N}{\log^{c+1} N} \right) = \tilde{O} \left( \log N \right)
\end{align*}
and
\begin{align*}
	R_{MR}(N) := \frac{T_{MR}(N)}{|\log \epsilon_{MR}(N)|} = \tilde{O}\left( \frac{\log^2 N}{2} \right) = \tilde{O}(\log^2 N )
\end{align*}

Since $c$ is a constant, $M(\log^{c+1}N) = \Theta(\log^c N \cdot M(\log N))$. Thus, 
\begin{align*}
	R_{AB}(N) = \tilde{O}\left( \log N \right)
\end{align*}
and 
\begin{align*}
	R_{MR}(N) = \tilde{O}(\log^2 N),
\end{align*}
which suggests that the accuracy guaranteed for \Algoref{algorithm:both} is stronger than that guaranteed for the Miller-Rabin algorithm relative to the algorithms' runtimes. \par

\section{Constructing a Polynomial that is Irreducible over $\mathbb{Z}/p\mathbb{Z}$} \label{section:polyconstruct}

Methods from \cite[Sections 2, 3, and 8]{LP} yield a deterministic algorithm which either correctly determines $N$ to be composite or constructs a polynomial of sufficiently high degree over $\mathbb{Z}/N\mathbb{Z}$ that is irreducible over $\mathbb{Z}/p\mathbb{Z}$ for each prime factor $p$ of $N$. This construction is developed through the language of pseudofields, (commutative and unital) rings that have an endomorphism that resembles a power of the Frobenius automorphism of finite fields. Furthermore, the algorithm combines both \cite[Algorithm 3.1]{LP} and \cite[Algorithm 8.3]{LP}. The former constructs a period system for $N$ and the latter uses the period system to construct pseudofields of small prime degree then takes the tensor product of these pseudofields. In doing so, the algorithm computes the desired polynomial. \par

We will show, in \Cref{proposition:irredreduction}, that any pseudofield that the algorithm of \cite{LP} constructs is a field when reduced modulo the prime factors of $N$. Since a pseudofield is isomorphic to a  residue ring of the form $(\ZNx)/(f(x))$, this means that the polynomial $f(x)$ that the algorithm computes is irreducible modulo the prime factors of $N$. Thus, $f(x)$ allows us to apply the strong accuracy bounds guaranteed by \Cref{proposition:accuracy2} and \Cref{proposition:strongaccuracy}. We will also produce, in \Cref{proposition:polyconstructruntime}, an asymptotic upper bound to the time it takes to construct such a pseudofield. \par

There are, however, a few restrictions in using these methods. One is that $N$ must be sufficiently large, in particular greater than an effectively computable constant which \cite{LP} calls $c_4$. Moreover, $\deg f$ will satisfy $\deg f \in [D,2D)$ and $N > \deg f$ where $D > (\log N)^{\frac{46}{25}}$. In particular, the algorithm constructs $f$ so that $\deg f$ is on the same order as $D$. \par

In \Cref{subsection:periodsystem}, \Cref{subsection:pseudofield}, and \Cref{subsection:gaussianperiods}, we summarize the details from \cite{LP} about period systems, pseudofields, and Gaussian periods that we will need. 

\subsection{Constructing Period Systems} \label{subsection:periodsystem}
\begin{definition}
	Let $N$ be an integer greater than $1$. A \textbf{period pair} for $N$ is a pair $(r,q)$ of integers such that
	\begin{itemize}
		\item $r$ is a prime number not dividing $N$,
		\item $q$ divides $r-1$ and $q > 1$,
		\item the multiplicative order of $N^{(r-1)/q}$ modulo $r$ equals $q$. 
	\end{itemize}
	Furthermore, a \textbf{period system} for $N$ is a finite set $\mathcal{P}$ of period pairs for $N$ such that 
	\begin{itemize}
		\item $\gcd(q,q') = 1$ whenever $(r,q), (r',q') \in \mathcal{P}, (r,q) \neq (r',q')$, 
	\end{itemize}
	and the \textbf{degree} of $\mathcal{P}$ is $\prod_{(r,q) \in \mathcal{P}} q$, denoted $\deg \mathcal{P}$. 
\end{definition}

The period system constructed by \cite[Algorithm 3.1]{LP} consists of period pairs $(r,q)$ in which $r$ and $q$ are bounded above with respect to the chosen integer $D$

\begin{proposition}	\label{proposition:periodsystem}
	There are effectively computable positive integers $c_4, c_5$ such that, for each integer $N > c_4$ and each integer $D > (\log N)^{\frac{46}{25}}$, there exists a period system $\mathcal{P}$ for $N$ consisting of pairs $(r,q)$ with
	$$
		r < D^{\frac{6}{11}}, \quad q < D^{\frac{3}{11}}, \quad q \text{ prime},
	$$
	and with $D \leq \deg \mathcal{P} < D + D^{1 - 1/c_5(\log \log D)^2} < 2D$. \par
	Moreover, \cite[Algorithm 3.1]{LP}, when given $N > 1$ and $D > 0$, computes a period system $\mathcal{P}$ for $N$ with these properties and satisfying $\deg \mathcal{P} \in [D,2D)$ if and only if such a period system exists, which is the case if $N > c_4$ and $D > (\log N)^{\frac{46}{25}}$. The runtime of \cite[Algorithm 3.1]{LP} is $\tilde{O}(D + D^{\frac{6}{11}} \log N)$. 
	\begin{proof}
		See \cite[Proposition 2.15]{LP} and \cite[Proposition 3.2]{LP}.
	\end{proof}
\end{proposition}

\subsection{Basic Properties of Pseudofields} \label{subsection:pseudofield}
For more on pseudofields, see \cite[Sections 2 and 5]{LP}.

\begin{definition} \label{definition:pseudofield}
A \textbf{pseudofield} is an ordered pair $(A,\alpha)$ where $A$ is a ring of characteristic $N$ and $\alpha \in A$ for which there are a positive integer $d$, called the \textbf{degree}, and a ring automorphism $\sigma$, which we call an \textbf{automorphism of the pseudofield}, such that 
\begin{itemize}
	\item $\# A \leq N^d$
	\item $\sigma \alpha = \alpha^N$ 
	\item $\sigma^d \alpha = \alpha$ 
	\item $\sigma^{d/l} \alpha - \alpha \in A^\times$ for all prime factors $l$ of $d$. 
\end{itemize}
\end{definition}

Moreover, a pseudofield of characteristic $N$ is of the form $(\ZNx)/(f(x))$ for some polynomial $f(x) \in \ZNx$.

\begin{lemma} \label{lemma:pseudofieldrep}
Let $(A,\alpha)$ be a pseudofield. There is a unique monic polynomial $f \in \ZNx$ such that there is an isomorphism $(\ZNx)/(f(x)) \simeq A$ that maps the residue of $x$ in $(\ZNx)/(f(x))$ to $\alpha$. The degree of $f$ is the degree of $A$.
\begin{proof}
	See \cite[Proposition 2.6]{LP} .
\end{proof}
\end{lemma}


\begin{lemma} \label{lemma:pseudofieldreduce}
	Let $(A,\alpha)$ be a pseudofield of degree $d$ and let $\sigma$ be the automorphism of $A$ as characterized in \Cref{definition:pseudofield}. For each prime factor $p$ of $N$ there exists a unique $i \in \mathbb{Z}/d\mathbb{Z}$ such that $\beta^p \equiv \sigma^i \beta \pmod{pA}$ for all $\beta \in A$.
	\begin{proof}
		See \cite[Proposition 5.4(f)]{LP}.
	\end{proof}
\end{lemma}

\subsection{Constructing Pseudofields with Gaussian Periods} \label{subsection:gaussianperiods}

 We reiterate details about Gaussian periods from \cite[Section 8]{LP}. Let $r$ be a prime number not dividing $N$ and let $\zeta_r$ be the reduction of $x$ in the residue ring $(\ZNx)/(\Phi_r)$ where $\Phi_r = \sum_{i=0}^{r-1} x^i$. In particular, $\zeta_r^r = 1$, $\zeta_r \neq 1$, and the elements $\zeta_r^i$, where $0 \leq i < r$, form a basis of $(\mathbb{Z}/N\mathbb{Z})[\zeta_r]$ over $\mathbb{Z}/N\mathbb{Z}$. \par

For each integer $a$ that is not divisible by $r$, the ring $(\mathbb{Z}/N\mathbb{Z})[\zeta_r]$ has an automorphism $\sigma_a$ which maps $\zeta_r$ to $\zeta_r^a$. The group of these $\sigma_a$ is denoted as $\Delta$ and $\Delta$ is isomorphic to $\mathbb{F}_r^\times$ and hence is a cyclic group of order $r-1$. Furthermore, the elements $\tau \zeta_r$ form a basis of $(\mathbb{Z}/N\mathbb{Z})[\zeta_r]$ over $\mathbb{Z}/N\mathbb{Z}$. \par

Let $q$ be a prime number dividing $r-1$ and let $\Delta^q$ denote the subgroup $\{\tau^q: \tau \in \Delta\}$ of $\Delta$. Let $\eta_{r,q} = \sum_{\rho \in \Delta^q} \rho \zeta_r$ and let
\begin{align*}
	f_{r,q}(x) = \prod_{\tau \Delta^q \in \Delta/\Delta^q} (x-\tau \eta_{r,q}),
\end{align*}
which is a monic polynomial of degree $q$. The polynomial $f_{r,q}$ specifies a pseudofield of characteristic $N$ and degree $q$, see \cite[Proposition 8.1(b)]{LP}. Moreover, $f_{r,q}$ is irreducible over $\mathbb{Z}/p\mathbb{Z}$ where $p$ is any prime factor of $N$.

\begin{proposition} \label{proposition:smallpseudofieldreduced}
	Let $p$ be a prime factor of a positive integer $N$, let $r$ be a prime number not dividing $N$ and let $q$ be a prime number dividing $r-1$. Let $(A,\alpha)$ be a pseudofield with automorphism $\sigma$ and characteristic polynomial $f_{r,q}$. By \Cref{lemma:pseudofieldreduce}, there exists a unique $i \in \mathbb{Z}/q\mathbb{Z}$ such that $\beta^p \equiv \sigma^i \beta \pmod{pA}$ for all $\beta \in A$. Then, $i \not\equiv 0 \pmod{q}$. 
	\begin{proof}
		Since the pseudofield $A$ is isomorphic to the residue ring $(\ZNx)/(f_{r,q})$, the reductions of the two rings modulo $p$ are isomorphic, i.e. $A/pA \simeq (\Zpx)/(f_{r,q})$. The automorphism $\sigma$ on $A$ reduces to an endomorphism on $A/pA$ because all elements of $pA$ are of the form $pa$ for some $a \in A$, and $\sigma(pA) = p\sigma(A)$. Thus, the Frobenius map $\beta \mapsto \beta^p$ is a ring endomorphism of $(\Zpx)/(f_{r,q})$ because the Frobenius map is the $i$-th power of this reduction of $\sigma$. \par
		Suppose for contradiction that $i \equiv 0 \pmod{q}$. In particular, $\beta^p \equiv \sigma^0 \beta \equiv \beta \pmod{pA}$ for all $\beta \in A$, i.e.~the Frobenius map is the identity on $A/pA$. Suppose that $f_{r,q}$ factorizes as
		\begin{align*}
			f_{r,q}= \prod_{i=1}^k f_i^{e_i}
		\end{align*}
		where the $f_i$ are distinct irreducible polynomials over $\Zpx$ and $e_i$ are positive integers. By the Chinese remainder theorem,
		\begin{align*}
			(\Zpx)/(f_{r,q}) \simeq \prod_{i=1}^k (\Zpx)/(f_i^{e_i}).
		\end{align*}
		Since the Frobenius map is identity map of $A/pA \simeq (\Zpx)/(f_{r,q})$, it is also the identity map of $(\Zpx)/(f_i)$ for each $i$. Since $f_i$ is irreducible, $(\Zpx)/(f_i)$ is a finite field, and the finite field of $p$ elements is the only finite field on which the Frobenius map is the identity. Thus, $f_i$ is of degree $d$ for every $i$, so the roots of $f_{r,q}$ are all elements of $\mathbb{Z}/p\mathbb{Z}$. Recall that 
		\begin{align*}
			\eta_{r,q} = \sum_{p \in \Delta^{q}} \rho \zeta_r = \sum_{i \in H} \zeta_r^i
		\end{align*}
		where $H$ is the index $q$ subgroup of $\mathbb{F}_r^\times$. Since $\zeta_r$ is defined by the relation $\sum_{i=0}^{r-1} \zeta_r^i$, the root $\eta_{r,q}$ of $f_{r,q}$ is not an element of $\mathbb{Z}/p\mathbb{Z}$, which is a contradiction. Hence, $i \not\equiv 0 \pmod{q}$.  
	\end{proof}
	
\end{proposition}

Furthermore, \cite[Algorithm 8.3]{LP} either correctly determines $N$ to be composite or constructs a pseudofield of larger degree by taking tensor products of pseudofields specified by polynomials of the form $f_{r,q}$ via \cite[Proposition 7.4]{LP}. 

\begin{proposition} \label{proposition:tensor}
	Let $(A_1,\alpha_1)$ and $(A_2,\alpha_2)$ be pseudofields such that $A_1$ and $A_2$ both have characteristic $N$ and suppose that their degrees $d_1$ and $d_2$ satisfy $d_1,d_2 > 1$ and $\gcd(d_1,d_2) = 1$. Let $\sigma_1$ and $\sigma_2$ respectively be the automorphisms of these pseudofields. The tensor product $(A_1 \otimes_{\mathbb{Z}/N\mathbb{Z}} A_2, \alpha_1 \otimes \alpha_2)$ is a pseudofield of characteristic $N$ and degree $d_1d_2$ with automorphism $\sigma_1 \otimes \sigma_2$. Furthermore, there is an algorithm that runs in $\tilde{O}(d_1d_2 \log N)$ time that either finds a prime factor of $N$ that is at most $d_1d_2$ or constructs the tensor product of $(A_1,\alpha_1)$ and $(A_2,\alpha_2)$.
	\begin{proof}
		See \cite[Proposition 7.1]{LP} and \cite[Proposition 7.4]{LP}.
	\end{proof}
\end{proposition}

\begin{proposition} \label{proposition:irredreduction}
	\cite[Algorithm 8.3]{LP} either correctly determines $N$ to be composite or constructs a pseudofield which is a field modulo every prime factor $p$ of $N$; that is, the algorithm produces a polynomial $f(x) \in \ZNx$ such that $f(x)$ is irreducible over $\mathbb{Z}/p\mathbb{Z}$.
	\begin{proof}
		If \cite[Algorithm 8.3]{LP} finds a prime factor of $N$, then the prime factor is either $N$ itself, in which case we conclude that $N$ is prime, or the prime factor is not $N$, in which case we conclude that $N$ is composite. Assume that the algorithm does not find a prime factor of $N$. \par
		\cite[Algorithm 8.3]{LP} constructs a pseudofield of characteristic $N$ of sufficiently large degree by first constructing pseudofields specified by some polynomials $f_{r_1,q_1},\ldots,f_{r_k,q_k}$ of degree $q_1,\ldots,q_k$ respectively, where the $q_i$'s are all distinct primes. The algorithm then finds the tensor product of these pseudofields, i.e.~finds the polynomial determining the tensor product. \par
		Say that $(A_1,\alpha_1),\ldots,(A_k,\alpha_k)$ are the pseudofields specified by $f_{r_1,q_1},\ldots,f_{r_k,q_k}$ with automorphisms $\sigma_1,\ldots,\sigma_k$ respectively. Let $A$ be the tensor product of these pseudofields, say that $A$ is determined by the polynomial $f$, let $\alpha = \alpha_1 \otimes \cdots \otimes \alpha_k$, let $\sigma = \sigma_1 \otimes \cdots \otimes \sigma_k$, and let $q = \prod_{j=1}^k q_j$ so that $A$ has degree $q$. \par
		Further let $p$ be any prime factor of $N$.	By \Cref{lemma:pseudofieldreduce}, for each $1 \leq j \leq k$, there is some unique $i_j \in \mathbb{Z}/q_j\mathbb{Z}$ such that $\alpha_j^p \equiv \sigma_j^{i_j} \alpha_j \pmod{pA_j}$. Moreover, $i_j \not\equiv 0 \pmod{q_j}$ by \Cref{proposition:smallpseudofieldreduced}. Via the Chinese remainder theorem, let $i \in \mathbb{Z}/q\mathbb{Z}$ be the unique residue satisfying $i \equiv i_j \pmod{q_j}$. In particular, $i$ is coprime to $q$ and
		\begin{align*}
			\sigma^i(\alpha) &= \sigma_1^i \alpha_1 \otimes \cdots \otimes \sigma_k^i \alpha_k \\
											 &\equiv \sigma_1^{i_1} \alpha_1 \otimes \cdots \otimes \sigma_k^{i_k} \alpha_k \pmod{pA}\\
											&\equiv \alpha_1^p \otimes \cdots \otimes \alpha_k^p \\
											&\equiv \alpha^p. 
		\end{align*}
		Since $\alpha$ generates $A$, the automorphism $\sigma^i$ of $A$ reduces to the Frobenius map on $A/pA$. 
		\par
		Let $f(x) = \prod_{l=1}^m f_l(x)^{e_l}$ be the factorization of $f(x)$ over $\mathbb{Z}/p\mathbb{Z}$ where $f_l(x)$ is irreducible. In particular, the following series of isomorphisms holds:
		\begin{align*}
			A/pA \simeq (\Zpx)[\alpha] \simeq (\Zpx)/f(x) \simeq \prod_{l=1}^m (\Zpx)/f_l(x)^{e_l}/
		\end{align*}
		Thus, the field $(\Zpx)/f_l(x)$ is isomorphic to a residue ring $F$ of $A$. Furthermore, the residue of $\alpha$ in $F$ generates $F$ over $\mathbb{Z}/p\mathbb{Z}$, so $\deg f_l$ divides every positive integer $a$ for which $\alpha = \alpha^{p^a}$ holds in $F$. Furthermore,
		\begin{align*}
			\alpha = \left( \sigma^q \right)^i \alpha =  \sigma^{qi} \alpha = \left( \sigma^i \right)^q (\alpha) = \alpha^{p^q}.
		\end{align*}
		holds in $F$, so $\deg f_l \mid q$. \par
		Suppose for contradiction that $\deg f_l \neq q$. Since $i$ is coprime to $q$, there is some integer $b$ such that $bi \equiv 1 \pmod{q}$. In particular, $b$ is coprime to $q$ and 
		\begin{align*}
			\alpha^{p^b} = \left(\sigma^i \right)^b \alpha = \sigma^{ib} \alpha = \sigma \alpha
		\end{align*}
		holds in $F$. Moreover, there is some prime factor $\mathfrak{l}$ of $q$ such that $\deg f_l \mid \frac{q}{\mathfrak{l}}$. Note that
		\begin{align*}
			\sigma^{q/\mathfrak{l}} \alpha = \alpha^{p^{bq/\mathfrak{l}}},
		\end{align*}
		so $\deg f_l \mid \frac{bq}{\mathfrak{l}}$. Since $b$ is coprime to $q$ and hence to $\deg f_l$, we have $\deg f_l \mid \frac{q}{\mathfrak{l}}$. Thus,
		\begin{align*}
			\sigma^{q/\mathfrak{l}} \alpha = \alpha
		\end{align*}
		in $F$. However, $\sigma^{q/\mathfrak{l}} \alpha - \alpha$ must be a unit of $A$ and hence in $F$, so the above equation is a contradiction. Hence, $\deg f_l = q$, so $f(x)$ is irreducible over $\mathbb{Z}/p\mathbb{Z}$ as desired. 
	\end{proof}
\end{proposition}

We now describe the runtime of \cite[Algorithm 8.3]{LP}.
\begin{proposition} \label{proposition:constructruntime}
	When given an integer $N > 1$, and a period system $\mathcal{P}$ satisfying $N > \deg \mathcal{P}$, \cite[Algorithm 8.3]{LP} runs in time
	\begin{align*}
		\tilde{O}\left( \left( \deg \mathcal{P} + \sum_{(r,q) \in \mathcal{P}} qr \right) \log N \right)
	\end{align*}
	and either correctly declares $N$ composite or constructs a pseudofield of characteristic $N$ and degree $\deg \mathcal{P}$. 
	\begin{proof}
		See \cite[Proposition 8.4]{LP}. 
	\end{proof}
\end{proposition}

\subsection{Runtime Bounds}
To reiterate, the desired polynomial, which was irreducible modulo each prime factor of $N$, was constructed by first constructing a period system $\mathcal{P}$ through \cite[Algorithm 3.1]{LP} and then running \cite[Algorithm 8.3]{LP} with $\mathcal{P}$ as an input. We give the runtime for constructing the polynomial.
\begin{proposition} \label{proposition:polyconstructruntime}
	There is an algorithm, when given positive integers $N$ and $D$ satisfying $N > \max(c_4,2D)$ and $D > (\log N)^{\frac{46}{25}}$ for an effectively computable constant $c_4$, that either correctly determines $N$ to be composite or produces a polynomial $f(x) \in (\ZNx)$ that is irreducible over $\mathbb{Z}/p\mathbb{Z}$ for every prime factor $p$ of $N$ and whose degree is in $[D,2D)$. Moreover, this algorithm runs in time
	\begin{align*}
		\tilde{O} \left( D \log N \right).
	\end{align*}
	\begin{proof}
		By \Cref{proposition:periodsystem}, \cite[Algorithm 3.1]{LP} computes a period system $\mathcal{P}$ such that $\deg \mathcal{P} \in [D,2D)$ and whose period pairs $(r,q)$ satisfy
		\begin{align*}
			r < D^{\frac{6}{11}}, \quad q < D^{\frac{3}{11}}, \quad q \text{ prime}.
		\end{align*}
		Using $\mathcal{P}$, \cite[Algorithm 8.3]{LP} produces a polynomial $f(x) \in \ZNx$ that is irreducible over $\mathbb{Z}/p\mathbb{Z}$ for each prime factor $p$ of $N$. By \Cref{proposition:periodsystem} and \Cref{proposition:constructruntime}, \cite[Algorithm 3.1]{LP} and \cite[Algorithm 8.3]{LP} run in times
		\begin{align*}
			\tilde{O}\left( D + D^{\frac{6}{11}} \log N \right) \text{ and } \tilde{O} \left( \left( \deg \mathcal{P} + \sum_{(r,q) \in \mathcal{P}} qr \right ) \log N \right)
		\end{align*}
		respectively. Note that $\mathcal{P} = \Theta(D)$ and that $qr = O\left( D^{\frac{9}{11}} \right)$. In particular,
		\begin{align*}
			\sum_{(r,q) \in \mathcal{P}} qr = O \left( D^{\frac{9}{11}} \# \mathcal{P} \right).
		\end{align*}
		Since $q$ is prime for each $(r,q) \in \mathcal{P}$, we can bound $\# \mathcal{P}$ as $O\left( \log D \right)$. Therefore, \cite[Algorithm 8.3]{LP} runs in time
		\begin{align*}
			\tilde{O}\left( \left( D  + D^{\frac{9}{11}} \log D \right) \log N \right)
		\end{align*}
		and hence in time
		\begin{align*}
			\tilde{O} \left( D \log N \right).
		\end{align*}
		Note that $D \log N$ asymptotically dominates the runtime of \cite[Algorithm 3.1]{LP}, so the combined algorithm runs in $\tilde{O}\left( D \log N \right)$ time as desired. 
	\end{proof}
\end{proposition}

\section{Conclusion}
We summarize the accuracy and runtime of the modified Agrawal-Biswas primality test that uses the polynomial constructed by the algorithm in \cite{LP}. 
\begin{theorem} \label{theorem:summary}
	Let $c > \frac{46}{25}$ and let $N > \max(1,c_4)$ be an odd positive integer where $c_4$ is a certain effectively computable constant. There is a Monte Carlo probabilistic algorithm, which indicates $N$ to be prime or composite, always correctly indicates $N$ to be prime if $N$ is prime, and falsely indicates $N$ to be prime when $N$ is composite with a probability at most
	\begin{align*}
		\frac{1}{2^{(r-1) \log^c N} \cdot N^{\log^c N - r}}.
	\end{align*}
	Moreover, the algorithm runs in time $\tilde{O}(\log^{c+2} N)$, assuming that each random bit can be generated in $O(1)$ time, and uses $O(\log^{c+1} N)$ random bits.
	\begin{proof}
	The Monte Carlo algorithm first uses \cite[Algorithm 3.1]{LP} and \cite[Algorithm 8.3]{LP} to either correctly determine $N$ to be composite or produce a polynomial $f(x) \in \ZNx$ that is irreducible over $\mathbb{Z}/p\mathbb{Z}$ for every prime factor $p$ of $N$. The algorithm then uses \Algoref{algorithm:both} to test whether or not $N$ is prime. Note that \Algoref{algorithm:both} indicates $N$ to be prime if $N$ is prime because it only runs the Miller-Rabin and modified Agrawal-Biswas algorithms. On the other hand, if $N$ is composite, then \Algoref{algorithm:both} falsely indicates $N$ to be prime with a probability at most
	\begin{align*}
		\frac{1}{2^{(r-1) \log^c N} \cdot N^{\log^c N - r}}.
	\end{align*}
	by \Cref{proposition:strongaccuracy}. \par
	Computing $f(x)$ takes $\tilde{O}(\log^{c+1}N)$ time by \Cref{proposition:polyconstructruntime}. Moreover, \Algoref{algorithm:both} runs in $\tilde{O}(\log^{c+2}N)$ time by \Cref{proposition:strongruntime} and uses $O(\log^{c+1} N )$ bits. Therefore, the Monte Carlo algorithm indeed takes $\tilde{O}(\log^{c+2}N)$ time by \Cref{proposition:strongruntime} and uses $O(\log^{c+1} N )$ bits. 
	\end{proof}
\end{theorem}


\begin{thebibliography}{1}


\bibitem{AB}
	M.~Agrawal and S.~Biswas, \emph{Primality and Identity Testing via Chinese Remaindering}, {\bf Journal of the ACM}, 50(4):429-443, 2003.
\bibitem{AGP}
	W. R.~Alford, A.~Granville, and C.~Pomerance, \emph{There are Infinitely Many Carmichael Numbers}, {\bf Annals of Mathematics}, 139(3):703-722, 1994.
\bibitem{AKS}
	M.~Agrawal, N.~Kayal, N.~Saxena, \emph{PRIMES is in P}, {\bf Annals of Mathematics}, 160(2):781-793, 2004.
\bibitem{AM}
	A.O.L.~Atkin and F.~Morain, \emph{Elliptic curves and primality proving}, {\bf Mathematics of Computation}, 61:29-68, 1993.
\bibitem{DLP}
	I.~Damg\aa rd, P.~Landrock, and C.~Pomerance, \emph{Average case error estimates for the strong
probable prime test}, {\bf Mathematics of Computation}, 61:177-194, 1993.
\bibitem{ModernComputerAlgebra}
	J.~von zur Gathen and J.~Gerhard, \emph{Modern Computer Algebra}. \newblock Cambridge University Press, 3rd~ed., 2013.
\bibitem{GK}
	S.~Goldwasser and J.~Killian, \emph{Almost all primes can be quickly certified}, Proceedings of the {\bf Eighteenth ACM Symposium on the Theory of Computing}, 316-329, 1986.
\bibitem{LP}
	H.~W.~Lenstra jr and C.~Pomerance, \emph{Primality testing with Gaussian periods}, preprint, 2017.
\bibitem{Miller}
G.L.~Miller, \emph{Riemann's Hypothesis and Tests for Primality}, {\bf Journal of Computer and System Sciences}, 13:300-317, 1976.
\bibitem{Monier}
L.~Monier, \emph{Evaluation and comparison of two efficient probabilistic primality testing
algorithms}, {\bf Theoretical Computer Science}, 12(1):97-108, 1980.
\bibitem{Rabin}
M.O.~Rabin, \emph{Probabilistic algorithm for testing primality}, {\bf Journal of Number Theory}, 12:128-138, 1980.
%
%
%
%
%
%
%
%
%
%
%
%
%
%
%
%
%
%
%
%
%
%
%
%
%
%
%
%

%
%
%
%
%
 %


\end{thebibliography}
\end{document}